\documentclass[reqno]{amsart}
\usepackage{amsfonts,amsmath,amssymb}

\newcommand{\R}{\mathbb{R}}

\newcommand{\p}{\partial}

\begin{document}

\title{Shape derivative of the first eigenvalue of the $1$-Laplacian}

\author{Nicolas Saintier}
\address{Departamento de Matem\'atica, FCEyN UBA (1428), Buenos Aires, Argentina}
\email{nsaintie@dm.uba.ar}


\date{June 6th, 2007}

\begin{abstract}{We compute the shape derivative of the functional $\Omega\to\lambda_{1,\Omega}$, where $\lambda_{1,\Omega}$ denotes the first eigenvalue of the $1$-Laplacian on $\Omega$. As an application, we find that the ball is critical among the volume-preserving deformations.}
\end{abstract}

\maketitle

Let $\Omega$ be a smooth bounded domain in $\mathbb{R}^n$, $n \ge 2$. The $1$-Laplacian on $\Omega$ is the formal operator
$$\Delta_1u = -\hbox{div}\left(\frac{\nabla u}{\vert\nabla u\vert}\right)$$
we get by a formal derivation of $F(u) = \int_\Omega\vert\nabla u\vert\,dx$, or by letting $p \to 1$ in the definition of the $p$-Laplacian, $p > 1$.  
By analogy with the definition of the first eigenvalue of the p-Laplacian on $\Omega$, we define the first eigenvalue $\lambda_{1,\Omega}$ of the $1$-Laplacian on $\Omega$ by the minimization problem
\begin{equation*}%\label{FirstDef}
 \lambda_{1,\Omega} = \inf_{\begin{cases} u\in \dot H_1^1(\Omega) \\ \int_{\Omega} |u|\, dx = 1 \end{cases}} 
                      \int_\Omega |\nabla u|\, dx\hskip.1cm ,
\end{equation*}
where $\dot H_1^1(\Omega)$ is the closure of $C^\infty_0(\Omega)$ in the Sobolev space $H_1^1(\Omega)$ of functions in $L^1(\Omega)$ with one derivative in $L^1$. 

\medskip

The purpose of this paper is the study of the dependence of $\lambda_{1,\Omega}$ under regular perturbations by diffeomorphisms of $\Omega$, i.e. we want to compute the first variation, the so-called shape derivative, of the functional $\Omega\to \lambda_{1,\Omega}$. General results about the stability of $\lambda_{1,\Omega}$ under perturbations of $\Omega$ have been obtained in \cite{HS}. In particular the authors of \cite{HS} found the shape derivative of $\Omega\to\lambda_{1,\Omega}$ in the case of regular perturbations by diffeomorphisms close to homotheties. We want to extend this  result to the case of a general perturbation by diffeomorphisms. 

Let us recall some known facts about $\lambda_{1,\Omega}$ (see e.g. \cite{HS,KawSchu}).
A natural space to study $\lambda_{1,\Omega}$ is the space $BV(\Omega)$ of functions of bounded variations (see, for instance, 
\cite{AFP,EvaGar,Giu,Ziemer}). By standard properties of $BV(\Omega)$, we can also define $\lambda_{1,\Omega}$ by 
\begin{equation}\label{SecondDef}
 \lambda_{1,\Omega} = \inf_{\begin{cases} u\in BV(\Omega) \\ \int_{\Omega} |u|\, dx = 1 \end{cases}} 
  \int_\Omega |\nabla u|+ \int_{\partial\Omega} |u|\,dH^{n-1},
\end{equation}
where $|\nabla u|$ is the total variation of the measure $\nabla u$, and $H^{n-1}$ denotes the $(n-1)$-dimensional Hausdorf measure. We refer to \cite{Dem3} for a detailed proof of this assertion.
Note here that if $u \in BV(\Omega)$ and $\overline{u}$ is the extension of $u$ by $0$ in ${\mathbb R}^n\backslash\overline{\Omega}$, then $\overline{u} \in BV({\mathbb R}^n)$ and
\begin{equation}\label{ExtEqt}
\int_{\R^n} |\nabla\overline{u}| = \int_\Omega |\nabla u| + \int_{\p\Omega} |u|\, dH^{n-1}.
\end{equation}
By lower semicontinuity of the total variation and compactness of the embedding $BV(\Omega) \hookrightarrow L^1(\Omega)$, it easily follows from (\ref{ExtEqt}) that the infimum in (\ref{SecondDef}) is attained by some nonnegative $u \in BV(\Omega)$. Then $u$ is a solution of the equation $\Delta_1u = \lambda_{1,\Omega}$ in the sense that there exists $\sigma \in L^\infty(\Omega,{\mathbb R}^n)$, $\|\sigma\|_\infty \le 1$, such that
\begin{equation}\label{vp11}
  \begin{cases}
 -\hbox{div}~ \sigma = \lambda_{1,\Omega} \hskip.1cm ,\\
   \sigma\nabla u = |\nabla u| \hskip.1cm\hbox{in}\hskip.1cm \Omega\hskip.1cm ,\hskip.1cm\hbox{and} \\
   (\sigma\vec{n})u = -u \hskip.1cm\hbox{on}\hskip.1cm \p\Omega\hskip.1cm ,
  \end{cases}
\end{equation}
where $\vec{n}$ is the unit outer normal to $\p\Omega$, and $ \sigma\nabla u $ is the distribution defined by integrating by parts $\int_\Omega(\sigma\nabla u)v\,dx$ when $v \in C^\infty_0(\Omega)$ and $\hbox{div}~\sigma$ makes sense (eee e.g. \cite{Dem3,ACC}). 
We then say that $u$ is an eigenfunction for $\lambda_{1,\Omega}$.

We can also express $\lambda_{1,\Omega}$ in a more geometric way as an isoperimetric type problem. 
We recall that a set $C\subset\R^n$ is said of finite perimeter if its characteristic function $\chi_C$ belongs to $BV(\R^n)$. We then define the perimeter $|\p C|$ of $C$ as $\int_{\R^n} |\nabla \chi_C|$. Using the coarea formula, we can rewrite (\ref{SecondDef}) as 
\begin{equation}\label{Def3}
 \lambda_{1,\Omega} = \inf_{C\subset\overline{\Omega},~\chi_C\in BV(\R^n)} \frac{|\p C|}{|C|}.
\end{equation} 
We refer e.g. to \cite{LI} for a proof of this assertion. 
This infimum is attained by some set of finite perimeter, e.g. by a level-set of an extremal for (\ref{SecondDef}), called an eigenset or also a Cheeger's set. Note that minimizers for $\lambda_{1,\Omega}$ touch the boundary $\partial\Omega$ since, if not, we may blow it up by a factor larger than one, which would decrease $\lambda_{1,\Omega}$. Uniqueness and nonuniqueness results of eigensets are in \cite{FriKaw,StrZie}. Concerning regularity, possible references are \cite{Alm,Gio,GonMasTam,LI,StrZie}.

\medskip

To study the variations of $\lambda_{1,\Omega}$ with respect to smooth variations of $\Omega$, we consider a smooth vector field $V:\R^n\to \R^n$ and, for small $t\in\R$, the diffeomorphisms $T_t$ defined by~
\begin{equation}\label{DiffeomDefEqt}
 T_t(x) = x + tV(x), 
\end{equation}
and eventually the perturbed domains $\Omega_t=T_t(\Omega)$. We want to compute the derivative at $t=0$ of the map $t\to\lambda_{1,\Omega_t}$. 

Shape analysis is the subject of an intense research activity. We refer for example to \cite{Henrot} for an introduction to this field. The shape derivative of the first eigenvalue $\lambda_{p,\Omega}$ of the $p$-Laplacian, $p>1$, has been computed in \cite{GarSab,Lam}:
$$ \frac{d}{dt} \lambda_{p,\Omega_t|t=0} = -(p-1)\int_{\p\Omega} \left|\frac{\p u_p}{\p\nu}\right|^p (V,\nu)\,dH^{n-1},$$
where $u_p$ is the unique positive normalized eigenfunctions for $\lambda_{p,\Omega}$ and $\nu$ is the unit normal vector to $\p\Omega$. What could be the shape derivative of the first eigenvalue of the $1$-Laplacian is thus not obvious from this formula. 

Our result is the following: 

\noindent{\bf Theorem.}
{\it Let $\Omega$ be a smooth bounded open subset of $\mathbb{R}^n$, $V:\R^n\to\R^n$ a smooth vector-field and $\Omega_t = T_t(\Omega)$, where the 
 $T_t$'s are the diffeomorphisms defined by (\ref{DiffeomDefEqt}). Then 
 $$ \lambda_{1,\Omega_t}\to \lambda_{1,\Omega} $$ 
 as $t\to 0$. Moreover, if we assume that there exists a unique nonnegative eigenfunction $u\in BV(\Omega)$ for $\lambda_{1,\Omega}$ such that
 $\int_\Omega |u|\,dx = 1$, then $u = |A|^{-1}\chi_A$ for some eigenset $A\subset\bar\Omega$, and the map $t\to\lambda_{1,\Omega_t}$ is
 differentiable at $t = 0$ with 
 \begin{equation}\label{der}
  \frac{d}{dt}\lambda_{1,\Omega_t|t = 0} = 
  \int_{\partial^* A} \left(\text{div}\,V-(\nu,DV\nu) + \lambda_{1,\Omega} (V,\nu)\right)\,\frac{dH^{n-1}}{|A|},
 \end{equation}
 where $\nu$ is the Radon-Nykodym derivative of $\nabla\chi_A$ with respect to $|\nabla\chi_A|$ (i.e. $\nabla\chi_A = \nu|\nabla\chi_A|$
 as measures), and $\p^*A$ denotes the reduced boundary of $A$, i.e. the part of the boundary of $A$ at which can be defined a notion of unit normal vector in a measure theoretic sense (see e.g. \cite{AFP,EvaGar,Giu,Ziemer}).
}
\medskip

Let us assume that $\Omega\subset\R^n$ is strictly convex. We then know from \cite{CCN,StrZie} that there exists a unique Cheeger set $A\subset\bar\Omega$. Since $A$ has minimum perimeter among all the subsets of $\bar\Omega$ of finite perimeter and of volume $|A|$, it follows from \cite{StrZie} that $\p A$ is $C^{1,1}$. Hence the unit exterior normal vector to $\p A$ is defined $H^{n-1}$ a.e. and Lipschitz. Its components are thus differentiable at $H^{n-1}$ almost every point of $\p A$. Note that this vector coincides with $-\nu$ $H^{n-1}$-a.e..
The mean curvature $H_{\p A}$ of $\p A$ is defined by $H_{\p A}=-\text{div}_{\p A}\nu$, where $\text{div}_{\p A}$ denotes the tangential derivative on $\p A$. We can now write that  
$$ \text{div}\,V-(\nu,DV\nu)=\text{div}_{\p A}\,V=\text{div}_g\,V_{\p A}-H_{\p A}(V,\nu),$$
where $V_{\p A}$ denotes the tangential part of $V$, and $\text{div}_g$ the divergence operator of the manifold $(\p A,g)$, $g$ being the metric on $\p A$ induced by the Euclidean metric (see e.g. \cite{Henrot}). We can thus rewrite ({\ref{der}) as 
\begin{equation*}
\begin{split}
  \frac{d}{dt}\lambda_{1,\Omega_t|t = 0} 
  & = \int_{\p A}  (\text{div}_g\,V_{\p A}-H_{\p A}(V,\nu) + \lambda_{1,\Omega} (V,\nu)) \,\frac{dH^{n-1}}{|A|} 
\end{split}                               
\end{equation*}
and thus 
\begin{equation}\label{der2}
 \frac{d}{dt}\lambda_{1,\Omega_t|t = 0}  = \int_{\p A}  (\lambda_{1,\Omega} - H_{\p A})(V,\nu) \,\frac{dH^{n-1}}{|A|}. 
\end{equation}

When $A=\bar\Omega$, a situation that happens when $\Omega\subset\R^n$ is smooth convex and its curvature is less that $|\p\Omega|/((n-1)|\Omega|)$ (see \cite{KawLac} when $n=2$, and \cite{ACC} for an arbitrary $n$), formula (\ref{der2}) writes as 
$$ \frac{d}{dt}\lambda_{1,\Omega_t|t = 0}  = \int_{\p \Omega}  (\lambda_{1,\Omega}-H_{\p\Omega})(V,\nu) \,\frac{dH^{n-1}}{|\Omega|},$$  
where $\nu$ is the inner unit normal to $\p\Omega$. In the particular case when $\Omega$ is a ball, $A=\bar\Omega$ and $H_{\p\Omega}$ is constant, so that if we consider a deformation that preserves the volume, i.e. a vector-field $V$ such that $\text{div}\,V=0$, we get 
$$ \frac{d}{dt}\lambda_{1,\Omega_t|t = 0}  = 0.$$
Hence a ball is critical for such deformations. 

\medskip

The following section is devoted to the proof of the theorem.

\section*{Proof of the theorem}

To simplify the notations, we let $\lambda=\lambda_{1,\Omega}$ and $\lambda_t=\lambda_{1,\Omega_t}$. 

Using the change of variable formula for functions of bounded variations \cite{Giu}, we can rewrite $\lambda_t$ as 
\begin{equation}\label{eq1}
 \lambda_t = \inf_{v\in BV(\Omega)}
    \frac{\displaystyle \int_{\bar\Omega} |D(x,t).\nu_v|C(x,t)|\nabla \bar v|}
         {\displaystyle \int_\Omega C(x,t) |v|\,dx}, 
\end{equation}
where $D(x,t) = (DT_t(x))^{-1}$, $C(x,t) = |\text{det}\,(D T_t(x))|$, and $\nu_v$ is the Radon-Nikodym derivative of $\nabla v$ with respect to 
$|\nabla v|$. Recall that $\bar v$ denotes the extension of $v$ to $\R^n$ by $0$ - see (\ref{ExtEqt}).
As $|\nu_v| = 1$ $|\nabla v|$ - a.e, 
\begin{equation}\label{vpp1_7}
  \lambda_t \le \inf_{v\in BV(\Omega)} \frac{\displaystyle\int_{\bar\Omega} |D(x,t)|C(x,t)|\nabla \bar v|}{\displaystyle \int_\Omega C(x,t) |v|\,dx}.
\end{equation}
Since $|D(x,t)|, C(x,t)\to 1$ as $t\to 0$ uniformly in $x\in\bar\Omega $, we deduce from (\ref{vpp1_7}) that 
\begin{equation}\label{vp3_3}
 \limsup_{t\to 0} \lambda_t\le \lambda.
\end{equation}
We let $u_t\in BV(\Omega_t)$ be a nonnegative eigenfunction for $\lambda_t$ normalized by $\int_{\Omega_t} u_t\,dx=1$, and 
$v_t = u_t\circ T_t\in BV(\Omega)$. Then $(\bar v_t)$ is bounded in $BV(\mathbb{R}^n)$. Indeed if we denote by $\bar\nu_t$ the Radon-Nikodym derivative of $\nabla u_t$ with respect to $|\nabla u_t|$, we have 
\begin{equation}\label{vp3_2}
\begin{split}
 \int_{\R^n} |\nabla \bar v_t|
 & =   \int_{\bar\Omega_t} \left|(DT_t^{-1})^{-1}\bar\nu_t\right|\left|\text{det}\,DT_t^{-1}\right| |\nabla\bar u_t| 
  \le   (1+o(1))\int_{\bar\Omega_t} |\nabla\bar u_t|  \\
 & =   (1+o(1))\lambda_t, 
\end{split}
\end{equation}
and 
\begin{eqnarray}\label{vp3_1}
 \int_\Omega v_t\,dx = \int_{\Omega_t} u_t |\text{det}\,DT_t^{-1}|\,dx  =  1 + o(1). 
\end{eqnarray}
We can thus assume that the $v_t$'s converge to some nonnegative $v\in L^1(D)$ in $L^1(D)$ and a.e. in $D$, where $D$ is a smooth bounded open subset of
$\R^n$ containing both $\Omega$ and the $\Omega_t$'s. We denote by $u$ the restriction of $v$ to $\Omega$. Then $u\in BV(\Omega)$, $u\ge0$, and in view of (\ref{vp3_1}), 
$$\int_\Omega u\,dx = 1. $$
The lower semi-continuity of the total variation and (\ref{vp3_2}) then give 
$$\lambda \le \int_{\bar\Omega}|\nabla \bar u| \le \int_{\R^n} |\nabla v| 
\le \liminf_{t\to 0} \int_{\R^n} |\nabla\bar v_t| \le \liminf_{t\to 0} \lambda_t. $$
We then deduce with (\ref{vp3_3}) that 
$$ \lambda_t\to \lambda = \int_{\bar\Omega} |\nabla \bar u|, $$ 
as $t\to 0$, and then that 
\begin{equation}\label{vp_diff_5}
 \int_{\bar\Omega} |\nabla \bar v_t| \to \lambda=\int_{\bar\Omega} |\nabla \bar u|, 
\end{equation}
as $t\to 0$. This proves the first part of the theorem. 

\medskip

Let us note for a future use that, as in \cite{HS}, it follow from (\ref{vp_diff_5}) that $|\nabla \bar v_t|\to |\nabla\bar u|$ weakly in the sense
that 
\begin{equation}\label{vp_diff_cv1}
\int_{\R^n} \phi |\nabla\bar v_t|\to\int_{\R^n} \phi |\nabla\bar u|
\end{equation}
for any $\phi\in C(\R^n)$ with compact support. 

\medskip 

We now prove the differentiability of the map $t\to\lambda_t$ and the formula (\ref{der}). We use $u$ as a test-function in (\ref{eq1}) to estimate $\lambda_t$, so that 
\begin{equation}\label{vp_diff_1}
 \lambda_t - \lambda
 \le  \frac{\displaystyle \int_{\bar\Omega} |D(x,t)\nu| C(x,t) |\nabla \bar u|}{\displaystyle \int_\Omega C(x,t)u\,dx} 
    - \lambda ,
\end{equation}
where $\nu$ is the Radon-Nikodym derivative of $\nabla u$ with respect to $|\nabla u|$. Since $|\nu| = 1$ $|\nabla u|$-a.e., we can assume that $|\nu|=1$ everywhere. Direct computations give that
\begin{equation}\label{dev1}
 C(x,t) = \text{det}\,(DT_t(x)) = 1 + t\text{div}\,V(x) + o(t),
\end{equation}
and 
\begin{equation}\label{dev2}
 |D(x,t)\nu| = 1 - (\nu;DV(x)\nu)t + o(t), 
\end{equation}
where the $o(t)$ is uniform in $x$. Thus (\ref{vp_diff_1}) becomes
\begin{equation*}
\begin{split}
\lambda_t - \lambda
 & \le \frac{\displaystyle \lambda + t\int_{\bar\Omega} \left(\text{div}\,V - (\nu,DV\nu)\right) |\nabla \bar u| + o(t)}
            {\displaystyle  1 + t\int_\Omega u\,\text{div}\,V\,dx + o(t)} - \lambda \\
 & = ~t\left(\int_{\bar\Omega} \text{div}\,V\,(|\nabla\bar u|-\lambda u\,dx) - \int_{\bar\Omega} (\nu,DV\nu) |\nabla \bar u| \right) + o(t).
\end{split}
\end{equation*}
Hence  
\begin{equation}\label{vp_diff_6}
 \limsup_{t\to 0^+}~ \frac{\lambda_t - \lambda}{t} \le  
\int_{\bar\Omega} \{(\text{div}\,V - (\nu,DV\nu))|\nabla\bar u| - \lambda u\,\text{div}\,V\,dx \},
\end{equation}
and 
\begin{equation}\label{vp_diff_7}
 \liminf_{t\to 0^-}~ \frac{\lambda_t - \lambda}{t} \ge  
\int_{\bar\Omega} \{(\text{div}\,V - (\nu,DV\nu))|\nabla\bar u| - \lambda u\,\text{div}\,V\,dx \}. 
\end{equation}

\medskip

It remains to prove the opposite inequalities. We use $\bar v_t$ as a test-function to estimate $\lambda$, so that 
$$ \lambda_t - \lambda = \int_{\bar\Omega_t} |\nabla\bar u_t| - \lambda
  \ge \int_{\bar\Omega} |D(x,t)\nu_t| C(x,t)\,|\nabla\bar v_t| 
      - \frac{\displaystyle \int_{\bar\Omega} |\nabla \bar v_t|}{\displaystyle \int_\Omega v_t\,dx}, $$
where $\nu_t$ denotes the Radon-Nykodym derivative of $\nabla \bar v_t$ with respect to $|\nabla \bar v_t|$. As previously, we can assume that 
$|\nu_t| = 1$ everywhere. In view of (\ref{dev1}) and (\ref{dev2}), we obtain 
\begin{equation}\label{vp_diff_8}
 \lambda_t - \lambda
  \ge \int_{\bar\Omega} |\nabla\bar v_t| + t\int\left(\text{div}\,V - (\nu_t,DV\nu_t)\right)|\nabla\bar v_t| 
     - \frac{\displaystyle \int_{\bar\Omega} |\nabla \bar v_t|}{\displaystyle \int_\Omega v_t\,dx} + o(t).
\end{equation}
Since $\text{div}\,V\in C(\bar\Omega)$, we get thanks to (\ref{vp_diff_cv1}) that 
$$ \int_{\bar\Omega} \text{div}\,V\,|\nabla\bar v_t| \to \int_{\bar\Omega} \text{div}\,V\,|\nabla\bar u|. $$
Independently, 
$$ \int_\Omega v_t\,dx = \int_\Omega u_t|\text{det}\,DT_t^{-1}|\,dx  $$
with 
$$ \text{det}\,DT_t^{-1} = 1 - t\,\text{div}\,V + o(t),$$
so that 
$$ \int_\Omega v_t\,dx = 1 -t\int_{\Omega_t} u_t\,\text{div}\,V\,dx + o(t) = 1 - t\int_\Omega u\,\text{div}\,V\,dx + o(t).  $$
Thus 
\begin{equation*}
\begin{split}
 \frac{\displaystyle \int_{\bar\Omega} |\nabla \bar v_t|}{\displaystyle \int_\Omega v_t\,dx}
 & = \int_{\bar\Omega} |\nabla \bar v_t| + t\int_{\bar\Omega} |\nabla \bar v_t|\int_\Omega u\,\text{div}\,V\,dx + o(t) \\
 & = \int_{\bar\Omega} |\nabla \bar v_t| + t\lambda\int_\Omega u\,\text{div}\,V\,dx + o(t),
\end{split}
\end{equation*}
where the last equality follows from (\ref{vp_diff_5}). Hence (\ref{vp_diff_8}) becomes 
\begin{equation*}\label{eq5}
 \lambda_t - \lambda \ge t\int_{\bar\Omega} \text{div}\,V(|\nabla\bar u|-\lambda u\,dx) 
                         - t\int_{\bar\Omega}(\nu_t,DV\nu_t)\,|\nabla\bar v_t| + o(t). 
\end{equation*}
Eventually, in view of (\ref{vp_diff_5}) and the weak  convergence of $\nabla \bar v_t$ to $\nabla \bar u$, which follows from the $L^1$ convergence of the $\bar v_t$ to $\bar u$, we can apply Reshetnyak' theorem \cite{AFP,Reshetnyak,LuMo} to get that 
 $$ \int_{\bar\Omega} g(x,\nu_t(x)) |\nabla \bar v_t| \to  \int_{\bar\Omega} g(x,\nu(x)) |\nabla \bar u| $$
for any continuous function $g:\bar\Omega\times S \to \R$, where $S$ denotes the unit sphere of $\R^n$. In particular, 
$$ \int_{\bar\Omega}(\nu_t,DV\nu_t)\,|\nabla\bar v_t| = \int_{\bar\Omega}(\nu,DV\nu)\,|\nabla\bar u| + o(1).$$
Hence 
$$ \lambda_t - \lambda \ge t\int_{\bar\Omega} \{(\text{div}\,V - (\nu,DV\nu))|\nabla\bar u| - \lambda u\,\text{div}\,V\,dx \} + o(t), $$
and thus 
\begin{equation*}
 \limsup_{t\to 0^+} \frac{\lambda_t - \lambda}{t} 
  \ge \int_{\bar\Omega} \{(\text{div}\,V - (\nu,DV\nu))|\nabla\bar u| - \lambda u\,\text{div}\,V\,dx \},
\end{equation*}
and 
\begin{equation*}
 \liminf_{t\to 0^-} \frac{\lambda_t - \lambda}{t}  
 \le \int_{\bar\Omega} \{(\text{div}\,V - (\nu,DV\nu))|\nabla\bar u| - \lambda u\,\text{div}\,V\,dx \}.
\end{equation*}
Since by assumption $u$ is the unique normalized eigenfunctions for $\lambda$, we deduce from these two inequalities and (\ref{vp_diff_6}), (\ref{vp_diff_7}) that the map $t\to\lambda_t$ is differentiable at $t=0$ with 
$$ \frac{d}{dt}\lambda_{t|t=0} = \int_{\bar\Omega} \{(\text{div}\,V - (\nu,DV\nu))|\nabla\bar u| - \lambda u\,\text{div}\,V\,dx \}.$$
Eventually, as there always exists an extremal for the problem (\ref{Def3}), there exists a set of finite perimeter $A\subset\bar\Omega$ such that 
$u=|A|^{-1}\chi_A$. It then follows from geometric measure theory that $|\nabla \bar u|=|A|^{-1}H^{n-1}_{\lfloor\p^*A}$ (see e.g. \cite{AFP,EvaGar,Giu,Ziemer}).
The previous formula becomes
$$ \frac{d}{dt}\lambda_{t|t = 0} = 
  \int_{\partial^* A} (\text{div}\,V-(\nu,DV\nu))\,\frac{dH^{n-1}}{|A|} - \lambda \int_A \text{div}\,V\,\frac{dx}{|A|}, $$
from which we deduce (\ref{der}) using Green' formula for sets of finite perimeter. This ends the proof of the theorem.

\vspace{1cm}

\noindent\textit{\textbf{Acknowledgments.}} The author aknowledges the support of the grants FONCYT PICT 03-13719 (Argentina).

\end{document}